\documentstyle[12pt]{article}
\begin{document}
\title{An Eulerian-Lagrangian Approach for Incompressible
Fluids: Local Theory}
\author{Peter Constantin 
\\Department  of Mathematics\\The University of Chicago }
\maketitle
\newtheorem{thm}{Theorem}
\newtheorem{prop}{Proposition}

\vspace{.5cm}

\section{Introduction}
The three dimensional Euler equations are evolution
equations for the three velocity components $u(x,t)$,
\begin{equation}
\frac{\partial u}{\partial t} + u\cdot\nabla u + \nabla p 
= 0,\label{e}
\end{equation}
coupled with a fourth equation, $\nabla\cdot u = 0$,
expressing incompressibility. In this Eulerian 
formulation the velocity $u$ and pressure $p$ are 
recorded at fixed locations $x\in{\mathbf R}^3$. The velocities and
pressure vanish at infinity or are periodic. 
The pressure is determined using 
incompressibility. The equation is conservative and
the total kinetic energy, $\int |u|^2dx$
is a constant of motion. 

The Euler equations can be studied in terms of
the  vorticity (\cite{maj}). The vorticity is a vector 
$\omega = \nabla\times u$
corresponding
to the anti-symmetric part
of the gradient matrix $\nabla u$. It obeys a quadratic equation, whose
nature is such that the magnitude of the
vorticity may increase in time. If the amplification
is not rapid enough then a well-known  
criterion (\cite{bkm}) guarantees that no blow up can occur: if 
$$
\int_0^T \sup_{x}|\omega (x,t)|dt <\infty
$$
and the initial data are smooth, then the solution is smooth on the time interval $[0,T]$. The vorticity equation can be interpreted as
the vanishing of a commutator
\begin{equation}
\left [D_t, \Omega\right ] = 0\label{com}
\end{equation}
where
$$
D_t = \frac{\partial}{\partial t} + u\cdot\nabla
$$
is the material derivative and
$$
\Omega  = \omega \cdot\nabla.
$$
The characteristics of the
first order differential operator $\Omega$ are vortex lines;
the characteristics of the material derivative $D_t$ are
Lagrangian particle paths. The Lagrangian variables are the
path maps $a\mapsto X(a,t)$. The connection between the Lagrangian
description and the Eulerian one is given by the relations
$$
u(x,t) = \frac{\partial X(a,t)}{\partial t}, \,\, x=X(a,t).
$$
In this paper we discuss a description of the Euler
equations as a system of three coupled active vector equations.
The description concerns Lagrangian quantities computed in Eulerian variables
and thus combines the physical significance of the
Lagrangian description with the analytical
advantages of the Eulerian description.
The description bears similarities to the Clebsch variable representation. The 
Clebsch variables are a pair of 
scalars, $\theta, \varphi$ that are constant on particle paths and can
be used to re-construct the velocity via
$$
u^i(x,t) = \theta(x,t)\frac{\partial \varphi (x,t)}{\partial x_i} -
\frac{\partial n(x,t)}{\partial x_i}.
$$
This interesting representation is somewhat
restrictive: not all solutions can be represented in this manner. That is
because the Clebsch variables impose special constraints 
on helicity.
Helicity is the scalar product of velocity and vorticity $h
=u\cdot\omega$. Although $h$ itself is not conserved on particle paths, 
the integrals
$$
\int_T h(x,t)dx = c
$$
are constants of motion, for any vortex tube $T$. 
A vortex tube $T$ is a time evolving
region in space (not necessarily simply connected)
whose boundary is at each point parallel to the vorticity,
$\omega\cdot \nu = 0$ where $\nu$ is the normal to
$\partial T$ at $x\in\partial T$. The constants $c$ reflect the degree of
topological complexity of the flow (\cite{moff}) and in general are 
non-trivial but they vanish identically for flows that admit a
Clebsch variables representation. Indeed, for such
flows the helicity is the divergence of a field that is parallel
to the vorticity $h=-\nabla \cdot
(n\omega )$. Topological properties of streamlines and vortex tubes
are relevant to hydrodynamic stability
(\cite{arkhe}) and turbulence (\cite{cho}, \cite{frisch}).
The description of the flow that allows for arbitrary
vortex structures is based on formula (\ref{NL}) 
(\cite{serr}, \cite{gold}, \cite{goldet}, \cite{hunt}) that was
used for numerical computations. Somewhat related Hamiltonian 
formulations have been 
introduced by several authors (\cite{kuz}, \cite{ose}, \cite{butt}, 
\cite{peg}).

\section{Eulerian-Lagrangian Description}
The Lagrangian formulation of the Euler equations describes the flow
in terms of a volume preserving diffeomorphism, the
map $a \mapsto X(a,t)$. The curve $t \mapsto X(a,t)$ is the Lagrangian path
at label $a$ and obeys Newton's law
\begin{equation}
\frac{\partial^2 X(a,t)}{\partial t^2} = F_X(a,t).\label{l}
\end{equation}
The incompressibility condition for the map is
\begin{equation}
det\left (\nabla_aX\right ) = 1.\label{inco}
\end{equation}
The initial condition sets the labels at the initial time:
$$
X(a,0)= a.
$$
The forces $F_X$ in (\ref{l}) are 
\begin{equation}
F_X(a,t) = -(\nabla_x p)(X(a,t)) = -\left[\left(\nabla_aX(a,t)\right)^*\right]^{-1}
(\nabla_a\tilde{p})(a,t) \label{F}
\end{equation}
with $\tilde{p}(a,t) = p(X(a,t))$ and where $p$ is the 
Eulerian pressure. The notation $M^*$ means the transposed of the matrix $M$,
$(M^*)^{-1}$ its inverse.
Multiplying (\ref{l}) by $(\nabla_aX)^*$ we obtain
\begin{equation}
\frac{\partial^2 X(a,t)}{\partial t^2}\left(\nabla_aX(a,t)\right)^* =
 - (\nabla_a\tilde{p})(a,t)\label{la}
\end{equation}
or, on components
\begin{equation}
\frac{\partial^2 X^j(a,t)}{\partial t^2}\frac{\partial X^j(a,t)}{\partial
  a_i} = - \frac{\partial\tilde{p}(a,t)}{\partial a_i}.\label{lai}
\end{equation}
Pulling out a time derivative in the left-hand side we obtain
\begin{equation}
\frac{\partial}{\partial t}\left [\frac{\partial X^j(a,t)}{\partial t}\frac{\partial X^j(a,t)}{\partial
  a_i}\right ] = -\frac{\partial\tilde{q}(a,t)}{\partial
  a_i}\label{nlai}
\end{equation}
where
\begin{equation}
\tilde{q}(a,t) = \tilde{p}(a,t) -\frac{1}{2}\left
|\frac{\partial X (a,t)}{\partial t}\right |^2\label{q}
\end{equation}
We integrate (\ref{nlai}) in time, fixing the label $a$:
\begin{equation}
\frac{\partial X^j(a,t)}{\partial t}\frac{\partial X^j(a,t)}{\partial
  a_i} = u_{(0)}^{i}(a) - \frac{\partial\tilde{n}(a,t)}{\partial
  a_i}\label{nl}
\end{equation}
where
\begin{equation}
\tilde{n}(a,t) = \int_0^t\tilde{q}(a,s)ds\label{n}
\end{equation}
and
\begin{equation}
u_{(0)}(a) = \frac{\partial X(a,0)}{\partial t}\label{u0}
\end{equation}
is the initial velocity. Note that $\tilde{n}$ has dimensions
of circulation or of kinematic viscosity (length squared per time).
The conservation of circulation 
$$
\oint_\gamma\frac{\partial X(\gamma ,t)}{\partial t}\cdot d\gamma
= \oint_\gamma\frac{\partial X(\gamma ,0)}{\partial t}\cdot d\gamma
$$
follows directly from the form (\ref{nl}).
Let us consider
\begin{equation}
A(x,t) = X^{-1}(x,t)\label{dA}
\end{equation}
the ``back-to-labels'' map, and note that it forms a vector of active
scalars (an active vector)
\begin{equation}
D_t A = \frac{\partial A}{\partial t} + u\cdot\nabla A = 0.\label{aA}
\end{equation}
Turning to (\ref{nl}), multiplying by
$\left[\left(\nabla_aX(a,t)\right)^*\right]^{-1}$ and reading at $a = A(x,t)$
we obtain the formula
\begin{equation}
u^i(x,t) = \left (u_{(0)}^j(A(x,t))\right )\frac{\partial
  A^j(x,t)}{\partial x_i} - \frac{\partial n(x,t)}{\partial x_i}\label{NL}
\end{equation}
where
\begin{equation}
n(x,t) = \tilde{n}(A(x,t))\label{N}
\end{equation}
The equation (\ref{NL}) shows that the general Eulerian velocity can
be written in a
form that generalizes the Clebsch variable representation:
\begin{equation}
u = (\nabla A)^*B - \nabla n
\label{nc}
\end{equation}
where $B = u_{(0)}(A(x,t))$ is also an active vector
\begin{equation}
D_t B= 0.\label{aB} 
\end{equation}
Conversely, and somewhat more generally,
if one is given a pair of active vectors $A = (A^1(x,t),\cdots ,
A^{M}(x,t))$ and $B = (B^{1}(x,t), \cdots , B^{M}(x,t)) $ of arbitrary
dimension $M$, such that the
active
vector equations (\ref{aA}) and (\ref{aB}) hold and if $u$ is given by
\begin{equation}
u(x,t) = \sum_{k=1}^M B^k(x,t)\nabla_xA^k(x,t) -\nabla_x n
\label{ncM}
\end{equation}
with some function $n$, then it 
follows that
$u$ solves the Euler equations
$$
\frac{\partial u}{\partial t} + u\cdot\nabla u + \nabla\pi = 0
$$
where
$$
\pi = D_t n + \frac{1}{2}|u|^2
$$
Indeed, the only thing one needs is the kinematic commutation relation
\begin{equation}
D_t \nabla_x f = \nabla_x D_t f -(\nabla_x u)^*\nabla_xf\label{comm}
\end{equation}
that holds for any scalar function $f$. The kinematic commutation relation (\ref{comm})
is a consequence of the chain rule, so it requires no assumption
other than smoothness. Differentiating (\ref{ncM}) and using the active vector
equations (\ref{aA}, \ref{aB}) it follows that
$$
D_t(u) = -\sum_{k=1}^M((\nabla_x u)^*\nabla_xA^k)B^k -\nabla_x(D_t n)
+ 
(\nabla_xu)^*\nabla n =
$$
$$
 -\nabla_x(D_tn) -(\nabla_x u)^*\left [\sum_{k=1}^M(\nabla_x A^k)B^k
- \nabla_x n\right ] = 
$$
$$
-\nabla_x(D_tn) -(\nabla_x u)^*u =  -\nabla_x(\pi).
$$

\section{The Active Vector Formulation}

The previous calculations can be summarized as follows:
A function $u(x,t)$ solves the incompressible Euler equations if and
only if it can be represented in the form $u = u_A$
with
\begin{equation}
u_A^i(x,t) = \phi^m\left (A(x,t)\right )\frac{\partial
  A^m(x,t)}{\partial x_i} - \frac{\partial n_A(x,t)}{\partial
  x_i}\label{ua}
\end{equation}
and
\begin{equation}
\nabla\cdot u_A = 0
\label{div}
\end{equation}
where $A(x,t)$ solves the active vector equation
\begin{equation}
\left (\partial_t + u_A\cdot\nabla\right )A = 0,\label{av}
\end{equation}
with initial data
$$
A(x,0) = x.
$$
The function $\phi$ represents the initial velocity and the
function $n_A(x,t)$ is determined up to additive constants by the
requirement of incompressibility, $\nabla\cdot u_A = 0$:
$$
\Delta n_A(x,t) = \frac{\partial}{\partial x_i}\left \{
\phi^m(A(x,t))\frac{\partial A^m(x,t)}{\partial x_i}\right \}.
$$
The periodic boundary conditions are 
\begin{equation}
A(x+Le_j, t) = A(x,t) + Le_j;\quad n_A(x + Le_j,t) = n_A(x,t)
\label{aper}
\end{equation}
with $e_j$ the standard basis in ${\mathbf{R}}^3$. In this case
\begin{equation}
\delta_A(x,t) = x - A(x,t),
\label{thetA}
\end{equation}
$n_A(x,t)$, and $u_A(x,t)$ are periodic functions in each spatial
direction. One may consider also the case of decay at infinity,
requiring that $\delta_A$, $u_A$ and $n_A$ vanish sufficiently fast
at infinity. 
The equation of state (\ref{ua}, \ref{div}) can be written as 
\begin{equation}
u_A = {\mathbf P}  \left\{\phi^m\left (A(\cdot,t)\right )
  \nabla A^m(\cdot ,t)\right \} = {\mathbf P}\left\{\left (\nabla A\right )^*\phi(A)\right\} 
  \label{uwP}
\end{equation}
where
\begin{equation}
{\mathbf P} = {\mathbf 1} - \nabla \Delta^{-1}\nabla \cdot\label{P}
\end{equation}
is the Leray-Hodge projector (with appropriate boundary conditions)
on divergence free functions. 

The Eulerian pressure is determined, up to additive constants by
$$
p(x,t) = \frac{\partial n_A(x,t)}{\partial t} + u_A(x,t)\cdot\nabla n_A(x,t)
+ \frac{1}{2}|u_A(x,t)|^2.
$$
The Jacobian obeys
$$
det\left (\nabla A (x,t)\right) = 1.
$$
The vorticity
$$
\omega_A(x,t) = \nabla \times u_A
$$
satisfies the Helmholtz equation 
\begin{equation}
D_t^A\omega_A = \omega_A\cdot \nabla u_A
\label{helmh}
\end{equation}
and is given by the Cauchy formula
\begin{equation}
\omega_A(x,t) = \left [\nabla A(x,t)\right ]^{-1}\zeta(A(x,t))
\label{oiaz}
\end{equation}
where $\zeta = \nabla\times\phi$ is the initial vorticity.

The advantage of an active vector formulation is that $A$ has
conserved distribution, that is, for any function $\Phi$
$$
\int\Phi(A(x,t))dx = const;
$$
in particular
$\|A(\cdot, t)\|_{L^{\infty}_{loc}(dx)}$ is constant in time.

\section{Local existence}
The proof of local existence of solutions to the Euler
equations in the active vector formulation is relatively simple
and the result can be stated economically.  The well-known local existence
results in Lagrangian (\cite{em}) and Eulerian (\cite{k}) variables 
require more derivatives and thus use more restrictive function spaces.

\begin{thm}
Let $\phi$ be a divergence free $C^{1,\mu}$ periodic vector
 valued function of three variables. There exists a time
interval $[0,T]$ and a unique $C([0,T]; C^{1, \mu})$
spatially periodic vector valued 
function $\delta (x,t)$ such that
$$
A(x,t) = x + \delta (x,t)
$$
solves the active vector formulation of the 
Euler equations, 
$$
\frac{\partial A}{\partial t} + u\cdot\nabla A = 0,
$$
$$
u = {\mathbf {P}}\left \{(\nabla A(x,t))^*\phi (A(x,t))\right \} 
$$
with initial datum $A(x,0) = x$.
\end{thm}
The same result holds if one replaces periodic boundary conditions
with decay
at infinity. Differentiating the 
active vector equation (\ref{av})
we obtain the equation obeyed by the gradients
\begin{equation}
D_t^A\left (\frac{\partial A^m}{\partial x_i}\right ) + \frac{\partial
  u_A^j}{\partial x_i}\frac{\partial A^m}{\partial x_j} =
  0.\label{grada}
\end{equation}
It is useful to denote
\begin{equation}
{\mathbf P}_{jl} = \delta_{jl} -
\partial_j\Delta^{-1}\partial_l\label{pcompo}
\end{equation}
the matrix elements of the Leray-Hodge operator.
Differentiating in the representation (\ref{uwP}) and
using the fundamental property
$$
{\mathbf P}_{jl}\frac{\partial f}{\partial x_l} = 0
$$
we obtain 
\begin{equation}
\frac{\partial u_A^j}{\partial x_i} = {\mathbf P}_{jl}\left (Det\left
[\zeta(A); \frac{\partial A}{\partial x_i}; \frac{\partial A}{\partial
  x_l}\right ]\right).\label{nablua}
\end{equation}
Recall that the function $\zeta$ is the curl of $\phi$. 
This relation shows that the
gradient of velocity can be expressed without use of second order derivatives
of $A$ and is the key to local existence: the equation 
(\ref{grada}) can be seen as a cubic quasi-local equation on characteristics. 
Let us make these ideas more precise. 
We will consider the periodic
case first. We write $C^{j,\mu}$,  $j = 0, 1$ to denote the 
H\"{o}lder spaces of real valued functions that
are defined for all $x\in {\mathbf R}^3$ and are periodic with period 
$L$ in each direction.
We denote by $\|f\|_{0,\mu}$ the $C^{0, \mu}$ norm:
\begin{equation}
\|f\|_{0,\mu} = \sup_x |f(x)| + \sup_{x\neq y}\left \{|f(x)-f(y)|\left
(\frac{L}{|x-y|}
\right )^{\mu} \right \}\label{omu}
\end{equation}
and by $\|f\|_{1,\mu}$ the  $C^{1, \mu}$ norm:
\begin{equation}
\|f\|_{1,\mu} =  \|f\|_{0,\mu} + L\|\nabla f\|_{0,\mu}\label{onemu}
\end{equation}
where the notation $|\cdots |$ refers to
modulus, Euclidean norm, and Euclidean norm for matrices, as appropriate.

We break the solution of the problem in two parts, the map
$\delta \to u$ and the map $u\to\delta $. 
We denote the first one $W$. 
\begin{equation}
W[\delta, \phi](x,t) =  {\mathbf {P}}\left \{({\mathbf {I}} +
\nabla\delta (x,t) )^*\phi
(x + \delta (x,t))\right \} \label{W}
\end{equation}
This map is linear in $\phi$ but nonlinear in $\delta$. 
\begin{prop}
The map $W[\delta, \phi]$
maps
$$
W: (C^{1,\mu})^3\times (C^{1,\mu})^3 \to (C^{1,\mu})^3
$$
continuously.
There exist constants $C$ depending on $\mu$ alone so that
$$
\|W[\delta, \phi]\|_{0,\mu} \le C \|\phi\|_{0,\mu}\left \{1 + \|\nabla
\delta \|_{0,\mu}\right \}^2
$$
and
$$
\|\nabla W[\delta, \phi]\|_{0,\mu} \le C \|\nabla\times\phi\|_{0,\mu}\left \{1 + \|\nabla
\delta \|_{0,\mu}\right \}^3.
$$
hold for any $\delta \in \left (C^{1,\mu}\right )^3$, $\phi\in \left (C^{1,\mu}\right )^3$.
\end{prop}

For the {\bf proof} we note that $W$ is made up from a number of operations. 
The first operation is the composition $\phi (x) \mapsto \phi (x+\delta (x))$. 
For a fixed $\delta \in (C^{1,\mu})^3$ the map  $x\mapsto x+\delta $
is Lipschitz. Composition with a Lipschitz change of variables maps
$C^{0,\mu}$ into itself continuously (we say that it is a continuous 
endomorphism).
The joint continuity of $[\phi, \delta ]\mapsto \phi(x+\delta)$
in $C^{1,\mu}$ follows naturally. The second operation is a sum of products of
functions (a matrix applied to a vector). This is a continuous operation
because the
H\"{o}lder  spaces $C^{j,\mu}$, $j= 0, 1$ we chose are Banach algebras.
The third and last operation is
the linear operator $\mathbf{P}$, which is bounded
in H\"{o}lder spaces. 
We need to consider also derivatives of $W$. We use the formula
(\ref{nablua}) and note that the expression for the
gradient is made of similar operations as above and apply the
same kind of reasoning. This finishes the proof.

Time does not play any
role in this proposition because the equation of state $(\delta, \phi)
\mapsto W[\delta, \phi]$ is time independent. The second half of the
procedure does depend on time. Let us denote by $\Theta$ the map that
associates
to two continuous paths $t\mapsto \delta(\cdot, t)$ and
$t\mapsto \phi (\cdot, t)$ a new path $t\mapsto \theta$;
the path
$t\mapsto \theta = \Theta[\delta, \phi]$ is obtained by solving the 
partial differential equation
\begin{equation}
\frac{\partial \theta}{\partial t} + u\cdot\nabla \theta + u = 0
\label{eqtheta}
\end{equation}
where
$$
u = W[\delta (\cdot, t), \phi (\cdot, t)],
$$
periodic boundary conditions are imposed on $\theta$ and zero initial data 
$$
\theta (x,0) = 0
$$
are required.
The Euler equation requires only the use of a time independent $\phi$, but
allowing time dependent $\phi$ is very useful: one can thus treat more
equations, in particular the Navier-Stokes equation.
Let us consider the space
$$
{\cal P}_T = C([0,T], (C^{1,\mu})^3)
$$
of continuous $(C^{1,\mu})^3$ -valued paths defined on a time interval
$[0,T]$, endowed with the natural norm
$$
\|\theta\|_{1,{\cal P}} = \sup_t\|\theta(\cdot, t)\|_{1,\mu}.
$$
We will consider also the weaker norm
$$
\|\theta\|_{0,{\cal P}} = \sup_t\|\theta(\cdot, t)\|_{0,\mu}.
$$
$\Theta$ is nonlinear in both arguments.  

\begin{prop}
The map $\Theta [\delta, \phi]$ maps
$$
\Theta : {\cal P}_T\times {\cal P}_T\to {\cal P}_T
$$
and is continuous when the topology of the source space ${\cal P}_T\times {\cal P}_T$ is the
natural product $C^{1,\mu}$ topology and the topology of the target space 
${\cal P}_T$ is the weaker $C^{0,\mu}$ topology. 
Moreover, there exists a constant $C$ depending on $\mu$ alone
so that
$$
\|\nabla \theta (\cdot, t)\|_{0,\mu} \le \left
(\int_0^t\|\nabla u(\cdot, s)\|_{0,\mu}ds\right )\left \{\exp\{C\int_0^t\|\nabla
u(\cdot, s)\|_{0,\mu}ds\}\right \}
$$
holds for each $t \le T$ with 
$u = W[\delta ,\phi]$ and $\theta  = \Theta[\delta, \phi]$.
\end{prop}
Proposition 2 states that the map $\Theta$ is bounded but not that it is 
continuous in the strong $C^{1,\mu}$ topology. The {\bf proof}  
follows naturally from the idea to use the classical 
method of characteristics and ODE Gronwall type arguments. Similar ideas are
needed below in the the slightly more difficult proof of 
Proposition 3 and we will sketch them there and therefore
 we leave the details of the proof of Proposition 2 to the interested reader.

In order to proceed let us take now a fixed
$\phi$, take a small number $\epsilon>0$ and associate to it the set
$$
{\cal I}\subset {\cal P}_T
$$
defined by
$$
{\cal I} = \{\delta (x,t); \delta (x,0) = 0,  \|\nabla \delta(\cdot,
t)\|_{0,\mu}\le \epsilon, \forall t \le T\}.
$$
Combining the bounds in the two previous propositions
one can choose, for fixed $\phi$, a $T$ small enough so that
$$
\delta \mapsto \Theta[\delta ,\phi] ={\cal S}[\delta]
$$
maps
$$
{\cal S}: {\cal I}\to {\cal I}.
$$
Inspecting the bounds it is clear that it is sufficient to require 
$$
T\|\nabla \times \phi\|_{0,\mu} \le c\epsilon
$$
with an appropriate $c$ depending on $\mu$ alone. Leaving $\phi$, $\epsilon$ and $T$
fixed as above, the map ${\cal S}$ is Lipschitz in the weaker norm
$C^{0,\mu}$:
\begin{prop} There exists a constant $C$, depending on $\mu$ alone,
such that, for every $\delta_1, \delta_2 \in {\cal I}$, the Lipschitz
bound
$$
\|{\cal S}[\delta_1 ] - {\cal S}[\delta_1 ]\|_{0,{\cal P}} \le
C \|\delta_1 -\delta_2 \|_{0, {\cal P}}
$$ 
holds.
\end{prop}
It is essential that $\delta_j\in {\cal I}$, so that they are smooth
and  their gradients are small, but nevertheless
this is a nontrivial statement. An inequality of the
type
$$
\|{\cal S}[\delta_1 ] - {\cal S}[\delta_1 ]\|_{0,{\cal P}} \le
C \|\delta_1 -\delta_2 \|_{1, {\cal P}}
$$ 
is easier to obtain, but loses one derivative. This kind of loss of
one derivative is a well-known difficulty in general compressible 
hyperbolic conservation laws. The situation is complicated in addition 
by the fact that the constitutive law $W$ depends on gradients. 
As we shall
see, incompressibility saves one derivative. The heart
of the matter is
\begin{prop}
Let $\phi\in (C^{1,\mu})^3$ be fixed. There exists a constant depending on $\mu$ alone
so that
$$
\|W[\delta_1, \phi] - W[\delta_2, \phi]\|_{0\mu} \le C \|\delta_1 -\delta_2\|_{0,\mu}\|\phi\|_{1,\mu}
$$
holds for any $\delta_j \in C^{1,\mu}$ with $\|\delta_j\|_{1,\mu} \le 1$.
\end{prop}
One could use the condition
$\delta_j \in C^{1,\mu}$ with $\|\delta_j\|_{1,\mu} \le M$ 
but then $C$ would depend on $M$ also.

\noindent {\bf Proof of Proposition 4.} Denoting
$$
u = W[\delta_1, \phi] - W[\delta_2, \phi],
$$
$$
\delta = \delta_1 -\delta_2,
$$
$$
\psi (x) = \frac{1}{2}\left (\phi(x+ \delta_1(x)) + \phi(x +\delta_2(x))\right ),
$$
$$
v(x) = \phi(x+ \delta_1(x)) - \phi(x +\delta_2(x)),
$$
$$
\gamma = \frac{1}{2}(\delta_1 + \delta_2)
$$
we write
$$
u = u_1 + u_2
$$
with 
$$
u_1 = {\mathbf P}\left \{ (\nabla\delta )^{*}\psi\right\}
$$
and
$$
u_2 = {\mathbf P}\left \{ ({\mathbf I} + \nabla\gamma )^{*}v\right\}
$$
Now the bound
$$
\|u_2\|_{0,\mu} \le C\|\delta \|_{0,\mu}\|\phi\|_{1,\mu}
$$
is obtained in the same way as the bound in Proposition 1. (Actually 
$\phi$ Lipschitz is enough here.)
The dangerous term is $u_1$ because it contains $\nabla\delta$. But here we
can ``integrate by parts'' and write
$$
u_1 = -{\mathbf P}\left\{(\nabla \psi)^{*}\delta \right\}
$$
because of incompressibility. The matrix $\nabla\psi$ is bounded in
$C^{0,\mu}$ and the bound follows again easily, as the bounds in Proposition 1.
This ends the proof of Proposition 4. We draw the attention to the fact that
the presence of the $*$ (transpose) operation is essential
for the ``integration by parts'' to be allowed.

\noindent Returning to the {\bf proof of
Proposition 3} we denote $\theta_j = {\cal S}\delta_j$,
$u_j = W(\delta_j, \phi)$, $u = u_1-u_2$, $\theta = \theta_1-\theta_2$
and write
$$  
\frac{\partial \theta}{\partial t} +
\frac{u_1+u_2}{2}\cdot\nabla\theta
+ u\cdot\nabla\left (\frac{\theta_1+\theta_2}{2}\right ) + u = 0
$$
We consider the characteristics $X(a,t)$ defined by
$$
\frac{dX}{dt} = \frac{u_1+u_2}{2}(X,t), \,\,\, X((a,0) =a
$$
and note that in view of Proposition 1 and the assumption 
$\delta_j\in {\cal I}$, the characteristics are well defined for $0\le t\le T$,
their inverse $A(x,t) = X^{-1}(x,t)$ (the ``back-to-labels'' map) is 
defined too. Moreover,
$$
\sup_{t, a}\left |\frac{\partial X}{\partial a}\right | \le C
$$
and
$$
\sup_{t,x}\left |\frac{\partial A}{\partial x} \right | \le C
$$
holds with a constant $C$ depending on $\mu$ alone. Consider now the function
$$
F(x,t) =  u\cdot\nabla\left (\frac{\theta_1+\theta_2}{2}\right ) + u.
$$
Solving by the method of characteristics we obtain
$$
\theta (x,t) = - \int_0^t F(X(A(x,t),s),s)ds.
$$
Using Proposition 4 in conjunction with the bounds in Propositions 1 and 2
we see that $F(x,t)$ is bounded (uniformly in time) in
$C^{0,\mu}$:
$$
\sup_t\|F(\cdot, t)\|_{0,\mu} \le C\|\phi\|_{1,\mu}\|\delta\|_{0,{\cal P}}
$$
Compositions with the uniformly Lipschitz $X$ and $A$ are harmless 
and we obtain the desired result
$$
\|\theta\|_{0,{\cal P}} \le C\|\delta\|_{0,{\cal P}}.
$$
This ends the proof of Proposition 3.
The {\bf proof } of Theorem 1 follows now using successive 
approximations. Starting  with a first guess $\delta_1\in {\cal I}$
we define inductively
$$
\delta_{n+1} = {\cal S}\delta_n \in{\cal I}.
$$
Proposition 3 implies that the sequence $\delta_n$ converges rapidly
in the $C^{0,\mu}$ topology to a limit $\delta$.  Because ${\cal I}$ is convex
it contains this weaker limit point, $\delta\in {\cal I}$. Because ${\cal S}$ 
has the weak Lipschitz property of Proposition 3 it follows that
${\cal S}\delta = \delta$. This actually means that $A=x+\delta(x,t)$ solves
the active vector formulation of the Euler equations and that $u= W[\delta, \phi]$ solves the usual Eulerian formulation.

Now let us consider the case of decay at infinity. This case is instructive
to look at this case because it illuminates
the difference between $\phi, u, W$
on the one hand and $x, \delta, \Theta$ on the other hand;
the function spaces
need to be modified in a natural fashion to accommodate this difference.
The issue of decay at infinity is both a physical one -- the total
kinetic energy must be defined, and a mathematical one -- ${\mathbf P}$
must be defined. But apart from this, the decay at infinity requirement does
not hinder the proof in any respect. 
\begin{thm} Let $\phi$ be a $C^{1,\mu}$ velocity that is square
integrable
$$
\int |\phi (x)|^2dx<\infty
$$
and whose curl is integrable to some power $1<q<\infty$,
$$
\int |\nabla\times\phi (x)|^qdx<\infty.
$$
Then for 
$\epsilon$ sufficiently small there exists
a time interval $[0,T]$ and a 
$C^{1,\mu}$ function $\delta (x,t)$ 
such that
$$
\sup_t\|\nabla \delta (\cdot, t)\|_{0,\mu} \le \epsilon
$$
and such that $x+ \delta(x,t)$
solves the active vector formulation of the Euler equation. The velocity
corresponding to this solution belongs to $C^{1,\mu}$, is square integrable
and the vorticity is integrable to power $q$.
\end{thm}
The proof follows along the same lines as above. 
Because $\phi$ enters linearly in the expression for $W$ and because we
control $\nabla\delta$ uniformly,  
issues of decay at infinity of  do not 
arise. In other words, the function space for velocities
does not need to be a Banach algebra, rather a module over the 
Banach algebra of the $\delta$ variables, which 
 need not decay at infinity. 

\section{The blow up issue}
Any solution of the Euler equation can be constructed using
a sequence of near identity transformations.
One starts out with
$$
\phi = u_0  
$$
and solves for an interval of time $0\le t\le t_1$ the active vector equation
$$
\frac{\partial A}{\partial t }+ u_A\cdot\nabla A= 0
$$
$$
u_A = {\mathbf P}\left ((\nabla A)^*u_{0}(A)\right )
$$
$$
A(x,0) = x.
$$
At time $t =t_1$ one resets:
$$
\phi = u_1 = u_A(\cdot , t_1)
$$
and solves the system above again, for a new time interval $t_1\le t\le t_2$
and so one continues the solution. The local existence result guarantees that
$$
(t_{n+1} - t_n )\|\nabla u_{n}\|_{0,\mu} \ge c>0
$$
and during this time the solution $A(x,t)$ remains close to the
identity in the sense that $\delta = A - x$ obeys
$$
\|\nabla \delta (\cdot ,t)\|_{0,\mu} \le \epsilon 
$$
with a prescribed $\epsilon <<1$. The formula $\ref{nablua}$ implies then that
$$
\|\nabla u_{n}\|_{0,\mu} \le K^n \|\nabla u_{0}\|_{0,\mu}
$$
with a fixed $K>1$. If the
inequalities above would be sharp then, of course, the time steps
would have to decrease exponentially and the procedure would diverge
in finite time. It is possible that for certain initial data the bounds may be 
overly pessimistic and
the solution may exist for a long time. But with the present knowledge, if one
desires long-lived solutions for arbitrary three dimensional data then
one needs to smooth either at the end of each step or during each time step. 
If one applies a smoothing procedure one evidently changes the problem
and one introduces an artificial dissipation. There are many ways one could
conceivable regularize the Euler equations. The physically correct 
energy  dissipating equation is the 
Navier-Stokes equation. Unfortunately it is not known in 
three dimensions if the Navier-Stokes equations have 
globally defined unique solutions that converge to solutions of the Euler
equations. Even in two dimensions, where the
existence of smooth solutions is known for both the Euler and Navier-Stokes
equations, the situation is not entirely trivial
(\cite{bc}, \cite{cw}). The two dimensional situation is characterized by
the absence of vortex stretching. In the case of
the three dimensional Euler equations the vorticity magnitude evolves 
according to the stretching
equation
\begin{equation}
D_t \left (|\omega|\right ) = \alpha |\omega|.\label{stretch}
\end{equation}
The stretching factor $\alpha $ is related to the vorticity
magnitude through a principal value singular integral (\cite{csiam}):
\begin{equation}
\alpha (x,t) = P.V. \int D\left(\hat{y}, \xi (x,t), \xi (x+y,t)\right )
|\omega (x+y,t)|\frac{dy}{|y|^3}.\label{alpha}   
\end{equation}
Here $\hat{y}$ is the unit vector in the direction of
$y$, $\xi (x,t) = \frac{\omega}{|\omega|}$ is the unit vector tangent to
the vortex line passing through $x$ at time $t$ and $D$ is a certain
geometric factor. The geometric factor is a smooth function of three
unit vectors, has zero average on the
unit sphere, $\int D dS(\hat{y}) = 0$ and vanishes pointwise
when $\xi (x,t) = \pm \xi(x+y,t)$.  
Because $\alpha $ has the same order of magnitude as $|\omega |$,
dimensional reasoning suggests blow up of the type one encounters
in the ordinary differential equation $\frac{dm}{dt} = m^2$,
$$
\sup_{x}|\omega (x,t)| \sim \frac{1}{T-t}.
$$
But if the vorticity direction $\xi$ is smooth then a geometric depletion of
$\alpha$  occurs; that means that  $\alpha$ is of the
order of magnitude of velocity times the magnitude of the
spatial gradient of $\xi$ (an inverse
length scale, assumed to be finite). The two dimensional Euler equations
correspond to the case $\xi = (0,0,1)$ and $\alpha = 0$ identically.
If
$$
\int_0^T\sup_{x}|\alpha (x,t)|dt <\infty
$$
then no blow up can occur. 
This idea of 
geometric depletion of nonlinearity
has been investigated theoretically and numerically for the Euler
equations and for a quasi-geostrophic active scalar equation 
(\cite{csiam}, \cite{cmt}, \cite{cfm}, \cite{oy}, \cite{cord}, \cite{cns}).
In the Eulerian-Lagrangian formulation of the Euler equation the
role played by smooth stratifications  
can be explained in the following manner. Consider functions 
$w = w_{\psi}$ of the form
\begin{equation}
w_{\psi}(x,t) = (\nabla A(x,t))^*\psi((A(x,t))\label{wA}
\end{equation}
associated to arbitrary vectors $\psi$. Alternately, one might
 consider solutions of
\begin{equation}
D_t^A w + \left (\nabla u_A\right )^* w = 0
\label{w}
\end{equation}
with initial data $\psi$.
A particular example is provided by choosing $\psi = \phi$, i.e. $w_{\phi} =
w_A$ 
\begin{equation}
w_A = (\nabla A)^*\phi(A)\label{Aw}
\end{equation}
which obeys
\begin{equation}
\nabla\times w_A = \nabla\times u_A= \omega_A.\label{curleq}
\end{equation}
Because the vorticity 
\begin{equation}
\omega_A = \nabla\times u_A\label{omegA}
\end{equation}
satisfies
\begin{equation}
\omega_A = (\nabla A)^{-1}\zeta (A)\label{omegaz}
\end{equation}
it follows that
\begin{equation}
\omega_A(x,t)\cdot w_{\psi}(x,t) = \zeta(A(x,t))\cdot
\psi(A(x,t)), \label{omegaw}
\end{equation}
holds for any $\psi$ or, in other words
\begin{equation}
D_t(\omega_A\cdot w) = 0\label{hel}
\end{equation}
holds for any solution $w(x,t)$ of (\ref{w}).
Global regularity of a solution of the Euler equations
would follow from  (\ref{omegaw}) if one could find
a {\it sufficient } family of vectors $\psi$. 
By a {\it sufficient} family for the 
initial velocity $\phi$ and the time interval
$[0,T]$ we mean a family of vectors $\psi$ such that 
there exists a non-negative
function $\gamma (t)$ with $\int_0^T\gamma dt< \infty$ such that
$$
|\omega_A(x,t)| \le \gamma (t)\sup_{\psi}|\omega_A(x,t)\cdot w(x,t)|
$$
holds for $0\le t\le T$. A sufficient family for all two
dimensional flows  is provided by just one $\psi$, 
$\psi = (0,0,1)$ with $\gamma = 1$. Generalizations would consist
of situations in which one could find sufficient families that depend on the
initial data and time and take locally the role played in 2D by
the vertical direction.

The blow up issue becomes, in terms of $A$, a question of formation
of infinite gradients in conserved quantities. This is similar to
the case of hyperbolic conservation laws but with the significant
difference that the underlying characteristic flow is volume-preserving:
$det(\nabla A) = 1$, the matrix $\nabla A$ is invertible and 
\begin{equation}
(\left (\nabla A(x,t)\right )^{-1})_{ij} = \frac{1}{2}\epsilon_{imn}Det
\left [e_j;
\frac{\partial A}{\partial x_m}; \frac{\partial A}{\partial x_n}
\right ]
\label{inv}
\end{equation}
holds, where $e_j =(\delta_{jk})$ is the canonical basis in ${\mathbf R}^3$.
Consider the Euler-Lagrange label differentiation
\begin{equation}
L^A_j = \frac{1}{2}\left  (\epsilon_{imn}\epsilon_{jkl}\frac{\partial
A_k}{\partial x_m}\frac{\partial A_l}{\partial x_n}\right
)\frac{\partial}{\partial x_i}\label{laj}
\end{equation}
From the commutation relation (\ref{comm}) and the A equation
(\ref{av}) it follows that
\begin{equation}
\left [ D_t^A, L_j^A \right ] = 0
\label{commu}
\end{equation}
holds for any $j= 1, 2, 3$. This commutation relation simply says that in
Lagrangian coordinates, time and label derivatives commute.
Note, from the formulas (\ref{omegaz}) and (\ref{inv}) that
\begin{equation}
\frac{1}{2}\epsilon_{pil}Det\left
[\zeta(A); \frac{\partial A}{\partial x_i}; \frac{\partial A}{\partial
  x_l}\right ] = \omega_A^p.\label{combo}
\end{equation}
It is clear now that $D_t^A$ commutes with $\Omega_A = \omega_A\cdot\nabla$
because it is  represented in terms of $L_j^A$:
\begin{equation}
\Omega _A = \zeta_j(A)L_j^A.\label{omegal}
\end{equation}
Observe that, 
in view of the definition of the operators $L_j^A$,
\begin{equation}
L_j^A = \left (\left (\nabla A(x,t)\right )^{-1}\right )_{kj}\frac{\partial}{\partial x_k}\label{laqj}
\end{equation}
it follows that
\begin{equation}
\left (\left (\nabla A(x,t)\right )^{-1}\right )_{ij} = L_j^A[x_i];\label{laji}
\end{equation}
on the other hand
\begin{equation}
D_t^A(x_i) = u_A^i
\label{obvi}
\end{equation}
holds, so from the commutation relation (\ref{commu})
we obtain
\begin{equation}
D_t^A \left (\left (\nabla A(x,t)\right )^{-1}\right )_{ij} = L_j^A(u_A^i).\label{Qui}
\end{equation}
This equation, which could have been derived also directly from
(\ref{grada}), implies the  vorticity equation because of (\ref{oiaz}):
\begin{equation}
D_t^A\omega_A = \zeta_j(A)L_j^A(u_A) = \Omega_A(u_A).\label{omegaL}
\end{equation}
Because of the result in (\cite{bkm}) and (\ref{oiaz}), it
is clear that the finiteness of 
$$
\int_0^T\|\left (\nabla A(\cdot,t)\right )^{-1}\|_{L^{\infty}(dx)}dt
$$
implies regularity. Or, using (\ref{inv}), we deduce that
the finiteness of
$$
\int_0^T\|\nabla A(\cdot,t)\|^2_{L^{\infty}(dx)}dt
$$
implies regularity. Let us introduce now the matrix
\begin{equation}
C_{ij}^A(x,t; z) = \left (\nabla A(x+z,t)\right )_{im} \left (\left (\nabla A(x,t)\right )^{-1}\right )_{mj}.\label{Cij}
\end{equation}
and call it
the Euler-Lagrange calibrator. The formula
\begin{equation}
C_{ij}^A(x,t; z) = L_j^{A(x,t)}\left (A_i(x+z,t)\right )\label{CL}
\end{equation}
shows that the calibrator measures the response of the
Eulerian translation to an infinitesimal Lagrangian translation. Note that
\begin{equation}
C_{ij}^A(x,t; 0) = \delta_{ij}.\label{ze}
\end{equation}
The calibrator is a quotient of gradients at different
locations and therefore locally spatially uniform, 
temporally arbitrary changes like dilations do not affect it.
The vorticity equation (\ref{omegaL}) can be expressed in terms of the
Euler-Lagrange calibrator (\cite{aim}):
\begin{equation}
D_t^A\omega_A^i = \left \{\frac{1}{4\pi}P.V. \int D( \zeta,
C^A\zeta, C_{., p}^A)\sigma_{il}(\hat{z})\frac{dz}{|z|^3}\right
\}\frac{\partial A_p(x,t)}{\partial x_l}\label{vorc}
\end{equation}  
Note that
$$
\omega_A^i\frac{\partial A_p(x,t)}{\partial
x_i}= \zeta_p(A(x,t))
$$
is bounded. It is therefore  natural to conjecture
that the
smoothness of $C^A$ prevents finite time blow up
for the Euler equations. This conjecture is true for the
quasi-geostrophic
active scalar. The interested reader is referred to (\cite{aim}) for details.

The blow up question for the Euler equations remains open. Numerical 
calculations provide insight and hints, but the answer will have to be
analytical. The considerations above point towards a possible incompressible
dispersive effect that hinders blow up: as the gradients of $A$ become large
the resulting rapid (\ref{NL}) and non-uniform (\ref{Qui}) motion disperses
the large gradients. This might cause instability of blow up or perhaps 
its suppression.

{\bf Acknowledgments.} I thank Diego Cordoba, Charles
Fefferman and Julian Hunt for helpful comments. 
This research was supported in part by NSF- DMS9802611. Partial support of 
AIM and the hospitality of the Princeton Mathematics Department are 
gratefully acknowledged.

\end{document}